\documentclass[12pt]{amsart}
\usepackage{stmaryrd,txfonts,fullpage,amsthm,amsmath,xcolor,picins,graphicx}
\usepackage{pdfpages}
\usepackage{hyperref}\hypersetup{colorlinks,
  linkcolor={green!50!black},
  citecolor={green!50!black},
  urlcolor=blue
}

\def\arXiv#1{{\href{http://front.math.ucdavis.edu/#1}{arXiv:#1}}}
\def\qed{{\hfill\text{$\Box$}}}

\def\bbC{{\mathbb C}}
\def\bbZ{{\mathbb Z}}

\begin{document}

\title[Unitarity Of Gassner]{A Note on the Unitarity Property of the Gassner Invariant}

\author{Dror Bar-Natan}
\address{
  Department of Mathematics\\
  University of Toronto\\
  Toronto Ontario M5S 2E4\\
  Canada
}
\email{drorbn@math.toronto.edu}
\urladdr{\url{http://www.math.toronto.edu/~drorbn}}

\date{\today; first edition: June 29, 2014.}

\thanks{This work was partially supported by NSERC grant RGPIN 262178. The full \TeX\ sources are at \url{http://drorbn.net/AcademicPensieve/2014-06/UnitarityOfGassner/}. Updated less often: \arXiv{1406.7632}.}

\maketitle

\begin{abstract}
We give a 3-page description of the Gassner invariant (or representation) of braids (or pure braids), along with a description and a proof of its unitarity property.
\end{abstract}

The unitarity of the Gassner representation~\cite{Gassner:OnBraidGroups} of the pure braid group was discussed by many authors (e.g.~\cite{Long:LinReps, Abdulrahim:Faithfulness, KirkLivingstonWang:Gassner}) and from several points of view, yet without exposing how utterly simple the formulas turn out to be\footnote{Partially this is because the formulas are simplest when extended a ``Gassner invariant'' defined on the full braid group, but then it is not a representation and it is not unitary. Yet it has an easy ``unitarity property''; see below.}. When the present author needed quick and easy formulas, he couldn't find them. This note is written in order to rectify this situation (but with no discussion of theory). I was heavily influenced by a similar discussion of the unitarity of the Burau representation in~\cite[Section~3.1.2]{KasselTuraev:BraidGroups}.

\parpic[r]{\parbox{1in}{
  $b_0=\sigma_1\sigma_3^{-1}\sigma_2$:
  \newline\includegraphics[width=1in]{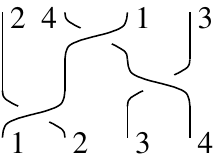}
}}
Let $n$ be a natural number. The braid group $B_n$ on $n$ strands is the group with generators $\sigma_i$, for $1\leq i\leq n-1$, and with relations $\sigma_i\sigma_j=\sigma_j\sigma_i$ when $|i-j|>1$ and $\sigma_i\sigma_{i+1}\sigma_i=\sigma_{i+1}\sigma_i\sigma_{i+1}$ when $1\leq i\leq n-2$. A standard way to depict braids, namely elements of $B_n$, appears on the right. Braids are made of strands that are indexed $1$ through $n$ at the bottom. The generator $\sigma_i$ denotes a positive crossing between the strand at position $\#i$ as counted just below the horizontal level of that crossing, and the strand just to its right. Note that with the strands indexed at the bottom, the two strands participating in a crossing corresponding to $\sigma_i$ may have arbitrary indices, depending on the permutation induced by the braids below the level of that crossing.

\parpic[r]{$\displaystyle
  U_{5;3}(t) = \begin{pmatrix}
    1 & 0 & 0 & 0 & 0 \\
    0 & 1 & 0 & 0 & 0 \\
    0 & 0 & 1-t & 1 & 0 \\
    0 & 0 & t & 0 & 0 \\
    0 & 0 & 0 & 0 & 1
  \end{pmatrix}
$}
Let $t$ be a formal variable and let $U_i(t)=U_{n;i}(t)$ denote the $n\times n$ identity matrix with its $2\times 2$ block at rows $i$ and $i+1$ and columns $i$ and $i+1$ replaced by $\begin{pmatrix} 1-t & 1 \\ t & 0 \end{pmatrix}$. Let $U^{-1}_i(t)$ be the inverse of $U_i(t)$; it is the $n\times n$ identity matrix with the block at $\{i,i+1\}\times\{i,i+1\}$ replaced by $\begin{pmatrix} 0 & \bar{t} \\ 1 & 1-\bar{t} \end{pmatrix}$, where $\bar{t}$ denotes $t^{-1}$.

\parpic[r]{$\displaystyle \Gamma(b)\coloneqq
  \prod_{\alpha=1}^k U_{i_\alpha}^{s_\alpha}(t_{j_\alpha})
$}
Let $b$ be a braid $b=\prod_{\alpha=1}^k \sigma_{i_\alpha}^{s_\alpha}$, where the $s_\alpha$ are signs and where products are taken from left to right. Let $j_\alpha$ be the index of the ``over'' strand at crossing $\#\alpha$ in $b$. The Gassner invariant $\Gamma(b)$ of $b$ is given by the formula on the right. It is a Laurent polynomial in $n$ formal variables $t_1,\ldots,t_n$, with coefficients in $\bbZ$.

\parpic[r]{\parbox{1.8in}{
  $\sigma_1\sigma_2\sigma_1$:\hspace{0.4in}$\sigma_2\sigma_1\sigma_2$:
  \newline\vspace{-3mm}
  \includegraphics[width=1.8in]{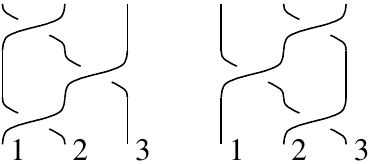}
}}
For example, $\Gamma(\sigma_1\sigma_2\sigma_1)=U_1(t_1)U_2(t_1)U_1(t_2)$ while $\Gamma(\sigma_2\sigma_1\sigma_2)=U_2(t_2)U_1(t_1)U_2(t_1)$. The equality of these two matrix products constitutes the bulk of the proof of the well-definedness of $\Gamma$, and the rest is even easier. The verification of this equality is a routine exercise in $3\times 3$ matrix multiplication. Impatient readers may find it in the {\sl Mathematica} notebook that accompanies this note,~\cite{Notebook}.

A second example is the braid $b_0$ of the first figure. Here and in~\cite{Notebook},
\[ \Gamma(b_0)=U_1(t_1)U_3^{-1}(t_4)U_2(t_1)
  = \begin{pmatrix}
    1-t_1 & 1-t_1 & 1 & 0 \\
    t_1 & 0 & 0 & 0 \\
    0 & 0 & 0 & \bar{t_4} \\
    0 & t_1 & 0 & 1-\bar{t_4}
  \end{pmatrix}
\]

\parpic[r]{$\displaystyle \Omega(\tau)\coloneqq
  \begin{pmatrix}
    (1-t_{\tau 1})^{-1} & 0 & \cdots & 0 \\
    1 & (1-t_{\tau 2})^{-1} & \cdots & 0 \\
    \vdots & \vdots & \ddots & \vdots \\
    1 & 1 & \ldots & (1-t_{\tau n})^{-1}
  \end{pmatrix}
$}
Given a permutation $\tau=[\tau 1,\ldots,\tau n]$ of $1,\ldots,n$, let $\Omega(\tau)$ be the triangular $n\times n$ matrix shown on the right (diagonal entries $(1-t_{\tau i})^{-1}$, $1$'s below the diagonal, $0$'s above). Let $\iota$ denote the identity permutation $[1,2,\ldots,n]$.

\noindent{\bf Theorem.} Let $b$ be a braid that induces a strand permutation $\tau=[\tau 1,\ldots,\tau n]$ (meaning, the strand indices that appear at the top of $b$ are $\tau 1,\tau 2,\ldots,\tau n$). Let $\gamma=\Gamma(b)$ be the Gassner invariant of $b$. Then $\gamma$ satisfies the ``unitarity property''
\begin{equation} \label{eq:unitarity}
  \Omega(\tau)\gamma^{-1}=\bar{\gamma}^T\Omega(\iota),
  \qquad\text{or equivalently,}\qquad
  \gamma^{-1}=\Omega(\tau)^{-1}\bar{\gamma}^T\Omega(\iota),
\end{equation}
where $\bar{\gamma}$ is $\gamma$ subject to the substitution $\forall i\, t_i\to \bar{t_i}\coloneqq t_i^{-1}$, and $\bar{\gamma}^T$ is the transpose matrix of $\bar{\gamma}$.

\noindent{\it Proof.} A direct and simple-minded computation proves Equation~\eqref{eq:unitarity} for $b=\sigma_i$ and for $b=\sigma_i^{-1}$, namely for $\gamma=U_i(t_i)$ and for $\gamma=U_i^{-1}(t_{i+1})$ (impatient readers see~\cite{Notebook}), and then, clearly, using the second form of Equation~\eqref{eq:unitarity}, the statement generalizes to products with all the intermediate $\Omega(\tau)^{-1}\Omega(\tau)$ pairs cancelling out nicely. \qed

If the Gassner invariant $\Gamma$ is restricted to pure braids, namely to braids that induce the identity permutation, it becomes multiplicative and then it can be called ``the Gassner representation'' (in general $\Gamma$ can be recast as a homomorphism into $M_{n\times n}(\bbZ[t_i,\bar{t_i}])\rtimes S_n$, where $S_n$ acts on matrices by permuting the variables $t_i$ appearing in their entries).

For pure braids $\Omega(\tau)=\Omega(\iota)\eqqcolon\Omega$ and hence by conjugating (in the $t_i\to1/t_i$ sense) and transposing Equation~\eqref{eq:unitarity} and replacing $\gamma$ by $\gamma^{-1}$, we find that the theorem also holds if $\Omega$ is replaced by $\bar{\Omega}^T$. Hence, extending the coefficients to $\bbC$, the theorem also holds if $\Omega$ is replaced by $\Psi\coloneqq i\Omega-i\bar{\Omega}^T$, which is formally Hermitian ($\bar{\Psi}^T=\Psi$).

If the $t_i$'s are specialized to complex numbers of unit norm then inversion is the same as complex conjugation. If also the $t_i$'s are sufficiently close to $1$ and have positive imaginary parts, then $\Psi$ is dominated by its main diagonal entries, which are real, positive, and large, and hence $\Psi$ is positive definite and genuinely Hermitian. Thus in that case, the Gassner representation is unitary in the standard sense of the word, relative to the inner product on $\bbC^n$ defined by $\Psi$.

We remark is that the Gassner representation easily extends to a representation of pure v/w-braids. See e.g.~\cite[Sections 2.1.2 and 2.2]{WKO1}, where the generators $\sigma_{ij}$ are described (they are {\em not} generators of the ordinary pure braid group). Simply set $\Gamma(\sigma_{ij})^{\pm 1}=U_{ij}^{\pm 1}$ where $U_{ij}$ is the $n\times n$ identity matrix with its $2\times 2$ block at rows $i$ and $j$ and columns $i$ and $j$ replaced by $\begin{pmatrix} 1 & 1-t_i \\ 0 & t_i \end{pmatrix}$. Yet on v/w-braids $\Gamma$ does not satisfy the unitarity property of this note and I'd be very surprised if it is at all unitary.

We also remark that there is an alternative form $\Gamma'$ for the Gassner representation of pure v/w-braids, defined by $\Gamma'(\sigma_{ij})^{\pm 1}=V_{ij}^{\pm 1}$ where $V_{ij}$ is the $n\times n$ identity matrix with its $2\times 2$ block at rows $i$ and $j$ and columns $i$ and $j$ replaced by $\begin{pmatrix} 1 & 1-t_j \\ 0 & t_i \end{pmatrix}$. Clearly, $U_{ij}$ and $V_{ij}$ are conjugate; $V_{ij}=D^{-1}U_{ij}D$, with $D$ the diagonal matrix whose $(i,i)$ entry is $1-t_i$ for every $i$. Hence on ordinary pure braids and for appropriate values of the $t_i$'s (as above), $\Gamma'$ is also unitary, relative to the Hermitian inner product defined by the matrix
\[ \Psi' \coloneqq \bar{D}^T\Psi D=i\bar{D}^T(\Omega-\bar{\Omega}^T)D \]
whose printed form is better avoided (yet it appears at the end of~\cite{Notebook}).

\AtEndDocument{\includepdf[pages={-}]{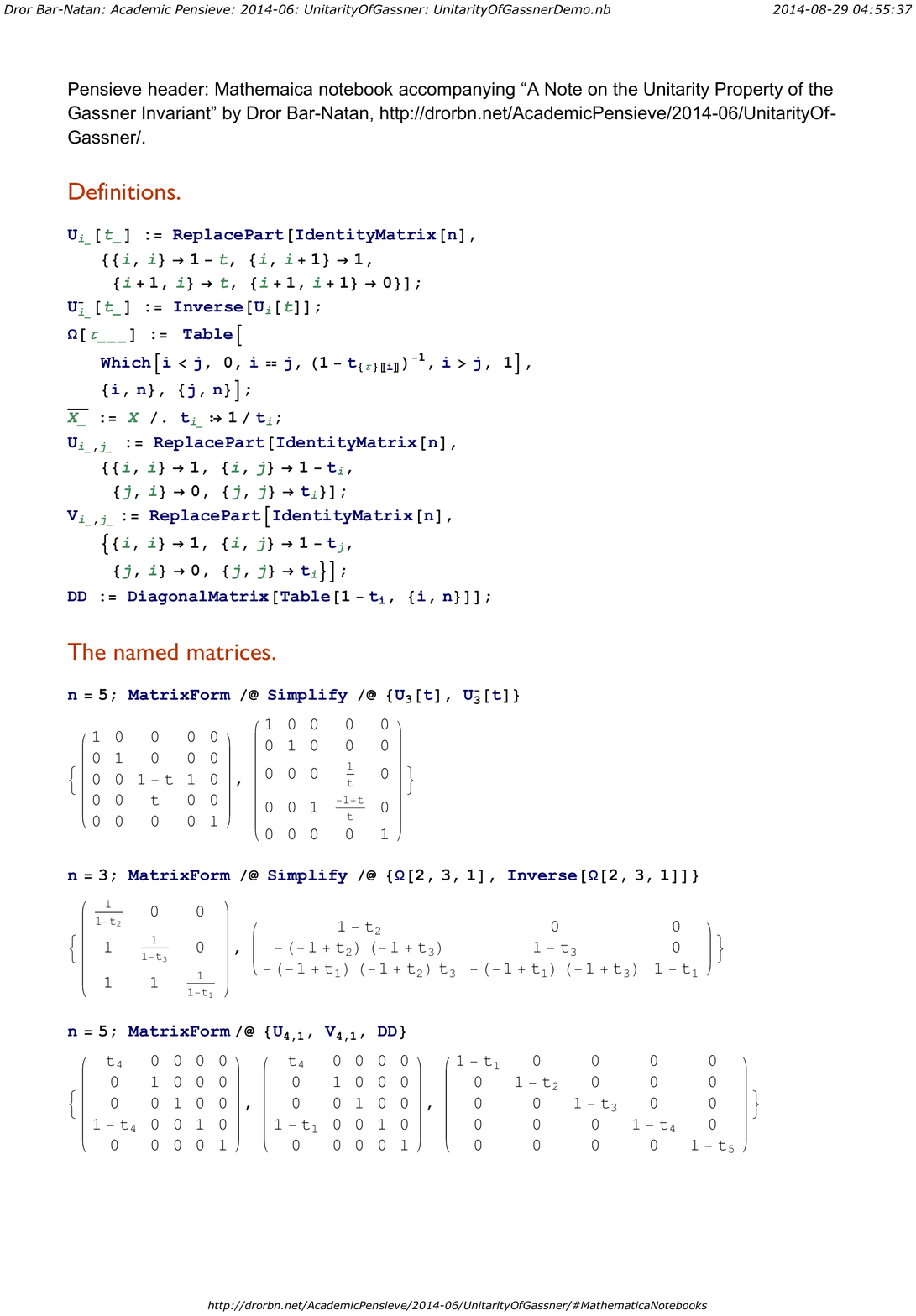}}


\begin{thebibliography}{KLW}

\bibitem[Ab]{Abdulrahim:Faithfulness} M.~N.~Abdulrahim,
  {\em A Faithfulness Criterion for the Gassner Representation of the Pure Braid Group,}
  Proceedings of the American Mathematical Society {\bf 125-5} (1997) 1249--1257.

\bibitem[BN]{Notebook} D.~Bar-Natan,
  {\tt UnitarityOfGassnerDemo.nb},
  a {\sl Mathematica} noteboook at \url{http://drorbn.net/AcademicPensieve/2014-06/UnitarityOfGassner/}.

\bibitem[BND]{WKO1} D.~Bar-Natan and Z.~Dancso,
  {\em Finite Type Invariants of w-Knotted Objects I: w-Knots and the Alexander Polynomial,}
  \url{http://drorbn.net/AcademicPensieve/Projects/WKO1/} and \arXiv{1405.1956}.

\bibitem[Ga]{Gassner:OnBraidGroups} B.~J.~Gassner,
  {\em On Braid Groups,}
  Ph.D. thesis, New York Univeristy, 1959.

\bibitem[KT]{KasselTuraev:BraidGroups} C.~Kassel and V.~Turaev,
  {\em Braid Groups,}
  Springer GTM {\bf 247}, 2008.

\bibitem[KLW]{KirkLivingstonWang:Gassner} P.~Kirk, C.~Livingston, and Z.~Wang,
  {\em The Gassner Representation for String Links,}
  Communications in Contemporary Mathematics {\bf 3-1} (2001) 87--136, \arXiv{math/9806035}.

\bibitem[Lo]{Long:LinReps} D.~D.~Long,
  {\em On the Linear Representation of Braid Groups,}
  Transactions of the American Mathematical Society {\bf 311-2} (1989) 535--560.

\end{thebibliography}
\end{document}